\documentclass[12pt]{amsart}
\usepackage{amsmath,amsthm,amscd,amsfonts,amssymb,graphicx,paralist,color}
\usepackage[a4paper,lmargin=2cm,rmargin=2cm,tmargin=4cm,bmargin=4cm]{geometry}
\usepackage{enumitem}
\usepackage[initials]{amsrefs}

\usepackage{xcolor}
\usepackage[hidelinks]{hyperref}
\hypersetup{colorlinks=true,  linkcolor=teal, citecolor=purple, urlcolor=orange}
\usepackage{rotating}

\usepackage[utf8]{inputenc}
\usepackage[T1]{fontenc}

\usepackage{tikz}

\usepackage{enumerate}

\theoremstyle{definition} 
\newtheorem{theorem}{Theorem}[section]
\newtheorem{lemma}[theorem]{Lemma} 
\newtheorem{corollary}[theorem]{Corollary} 
\newtheorem{proposition}[theorem]{Proposition}

\newtheorem{definition}[theorem]{Definition}

\newtheorem{remark}[theorem]{Remark}

\newcommand{\card}{\mathrm{card}}

\newcommand{\co}{\mathfrak{c}}

\newcommand{\nn}[1]{{\left\vert\kern-0.25ex\left\vert\kern-0.25ex\left\vert #1 
		\right\vert\kern-0.25ex\right\vert\kern-0.25ex\right\vert}}
\renewcommand{\geq}{\geqslant}
\renewcommand{\leq}{\leqslant}

\newcommand{\eps}{\varepsilon}

\newcounter{smallromans}

	{\end{list}}

\makeatletter
\renewcommand{\tocsection}[3]{%
	\indentlabel{\@ifnotempty{#2}{\bfseries\ignorespaces#1 #2\quad}}\bfseries#3}
\renewcommand{\tocsubsection}[3]{%
	\indentlabel{\@ifnotempty{#2}{\ignorespaces#1 #2\quad}}#3}
\newcommand\@dotsep{4.5}
\def\@tocline#1#2#3#4#5#6#7{\relax
	\ifnum #1>\c@tocdepth % then omit
	\else
	\par \addpenalty\@secpenalty\addvspace{#2}%
	\begingroup \hyphenpenalty\@M
	\@ifempty{#4}{%
		\@tempdima\csname r@tocindent\number#1\endcsname\relax
	}{%
		\@tempdima#4\relax
	}%
	\parindent\z@ \leftskip#3\relax \advance\leftskip\@tempdima\relax
	\rightskip\@pnumwidth plus1em \parfillskip-\@pnumwidth
	#5\leavevmode\hskip-\@tempdima{#6}\nobreak
	\leaders\hbox{$\m@th\mkern \@dotsep mu\hbox{.}\mkern \@dotsep mu$}\hfill
	\nobreak
	\hbox to\@pnumwidth{\@tocpagenum{\ifnum#1=1\bfseries\fi#7}}\par% <-- \bfseries for \section page
	\nobreak
	\endgroup
	\fi}
\AtBeginDocument{%
	\expandafter\renewcommand\csname r@tocindent0\endcsname{0pt}
}
\def\l@subsection{\@tocline{2}{0pt}{2.5pc}{5pc}{}}
\makeatother

\usepackage{todonotes}

\begin{document}

	\title[Algebraic genericity of certain families of nets in Functional Analysis]{Algebraic genericity of certain families of nets in Functional Analysis}
	%\today
	
	\author[Dantas]{Sheldon Dantas}
	\address[Dantas]{Departamento de Análisis Matemático, Facultad de Ciencias Matemáticas, Universidad de Valencia, Doctor Moliner 50, 46100 Burjasot (Valencia), Spain. \newline
		\href{http://orcid.org/0000-0001-8117-3760}{ORCID: \texttt{0000-0001-8117-3760} } }
	\email{\texttt{sheldon.dantas@uv.es}}
	
	\author[Rodríguez-Vidanes]{Daniel L. Rodríguez-Vidanes}
	\address[Rodríguez-Vidanes]{Instituto de Matemática Interdisciplinar (IMI), Departamento de Análisis Matemático y Matemática
		Aplicada, Facultad de Ciencias Matemáticas, Plaza de Ciencias 3, Universidad Complutense de Madrid,
		Madrid, 28040, Spain \newline
		\href{http://orcid.org/0000-0000-0000-0000}{ORCID: \texttt{0000-0003-2240-2855} }}
	\email{\texttt{dl.rodriguez.vidanes@ucm.es}}

	\begin{abstract} 
		In Functional Analysis, certain conclusions apply to sequences, but they cannot be carried over when we consider nets. In fact, some nets, including sequences, can behave unexpectedly. In this paper we are interested in exploring the prevalence of these unusual nets in terms of linearity. Each problem is approached with different methods, which have their own interest. 
		As our results are presented in the contexts of topological vector spaces and normed spaces, they generalize or improve a few ones in the literature.
		We study lineability properties of families of (1) nets that are weakly convergent and unbounded, (2) nets that fail the Banach-Steinhaus theorem, (3) nets indexed by a regular cardinal $\kappa$ that are weakly dense and norm-unbounded, and finally (4) convergent series which have associated nets that are divergent.
	\end{abstract}
	
	
	\maketitle
	
	\tableofcontents
	
	\thispagestyle{plain}
	
	\section{Introduction}
	
	The present paper is about lineability: we say that a subset $M$ of a vector space $X$ is {\it lineable} (respectively, $\kappa$-{\it lineable}, for a cardinal $\kappa$) if $M \cup \{0\}$ contains a vector space of infinite-dimension (respectively, of dimension $\kappa$). 
	These sort of problems have been studied intensively during the past few years since the term lineability was coined by V. I. Gurariy in the early 2000s and since then have appeared in the literature in different areas of Mathematics such as Functional Analysis, Real Analysis, Complex Analysis, Set Theory, Dynamical Systems, among others (we send the reader to \cites{ABPS, BPS, FGK, GQ, LRS, LM, Rm, SS} and the references therein).
	
	Our main interest here is to provide a contribution about the study of ``pathological'' nets %which (do not) satisfy relevant results 
	in Functional Analysis in the sense of lineability. 
	%More precisely, we invoke three (well-known) results from this area that hold true for sequences but when one tries to extend them to the context of nets difficulties arise immediately. 
	Some of these families of nets that are studied in this work arise from two well-known results that hold true for sequences.
	They are the following ones: (a) every weakly convergent sequence is bounded and (b) the Banach-Steinhaus theorem, which states that every sequence of bounded linear operators that converges pointwise to a bounded linear operator is uniformly bounded. On the other hand, for the upcoming results it is also important to recall also the following.
	Some authors have thoroughly worked on the problem about finding conditions that $(\|x_n\|)_{n \in \mathbb{N}}$ has to satisfy in order that the set $\{x_n: n \in \mathbb{N}\}$ is weakly closed (see, for instance, \cite{AGM, B, Kadets, KLO}). In particular, it is known that for every separable Banach space $X$, there exists a sequence $(x_n)_{n \in \mathbb N} \subseteq X$ such that $\|x_n\| \longrightarrow \infty$ as $n \rightarrow \infty$ and the set $\{x_n: n \in \mathbb{N} \}$ is weakly dense in $X$ (see, for instance, \cite[Corollary 5]{AGM} and also \cite[Section 2]{KLO} for a refinament).

	\vspace{0.3cm}
	\subsection{Preliminaries and notation} In what follows, we present briefly the notation we will be using throughout the paper and then we will describe our main results. All the spaces that we work with here are considered to be nonzero. We will be using basic concepts and notations from Set Theory found, for instance, in \cite{C, Je}. 
	Ordinal numbers will be identified with the set of their predecessors and cardinal numbers with the initial ordinals. 
	Given a set $A$, the cardinality of $A$ will be denoted by $\card(A)$. We denote by $\aleph_0$, $\aleph_1$ and $\co$ the first infinite cardinal, the second infinite cardinal and the cardinality of continuum, respectively.
	If $\alpha$ is an ordinal number, then its cofinality is denoted by $\text{cof}(\alpha)$.
	We say that a cardinal number $\kappa$ is regular if $\text{cof}(\kappa)=\kappa$.
	
	A set $\mathcal A$ is a directed set (also known as an index set) if $\mathcal A$ is a nonempty set that is endowed with a preorder $\leq$ (a reflexive and transitive relation) such that every pair of elements of $\mathcal A$ has an upper bound.
	A net in a set $X$ is a function from a directed set $\mathcal A$ to $X$ which will be denoted by $(x_a)_{a\in \mathcal A}$.
	We denote the set of nets in $X$ indexed by $\mathcal A$ as $X^{\mathcal A}$.
	%(observe that $\mathbb N$ with the standard order is a directed set and $(x_n)_{n\in \mathbb N}$ is the usual sequence).
	
	Given  a topological space $X$ and $(x_a)_{a\in \mathcal A}$ a net in $X$, we say that $(x_a)_{a\in \mathcal A}$ converges to $x\in X$ if for every neighborhood $U^x$ of $x$, there exists an element $a_0 \in \mathcal A$ such that $x_a \in U^x$ for every $a\geq a_0$.
	Recall that if $X$ is a topological vector space, then a net $(x_a)_{a\in \mathcal A}$ in $X$ weakly converges to $x\in X$ (denoted by $x_a \stackrel{w}{\longrightarrow} x$) if and only if $(x^*(x_a))_{a\in \mathcal A}$ converges to $x^*(x)$ for every $x^*\in X^*$.
	
	If a directed set $\mathcal A$ is in particular an ordinal number $\alpha$, then we have the so-called $\alpha$-sequences instead of nets defined in a set.
	It is known that the convergence of $\alpha$-sequences in a topological space can be reduced to the convergence of $\text{cof}(\alpha)$-sequences (see, for instance, \cite[Propositions~3.1 and~3.2]{R}).
	Therefore, we simply consider $\kappa$-sequences, where $\kappa$ is a regular cardinal number.
	This notion of $\kappa$-sequence was introduced in 1907 by J. Mollerup \cite{Mo} and has been studied throughout the 20th and 21st centuries by many mathematicians in several contexts (see \cite{Ko,Si,Na,R} and the references therein).
	Given a $\kappa$-sequence $(x_\alpha)_{\alpha < \kappa}$ in $X$, we say that $(x_{\beta_\alpha})_{\alpha < \kappa}$ is a $\kappa$-subsequence of $(x_\alpha)_{\alpha < \kappa}$ if there exists an increasing injective function $\varphi : \kappa \to \kappa$ such that $x_{\beta_\alpha} = x_{\varphi(\alpha)}$ for every $\alpha < \kappa$.
	
	As stated earlier, in this work we are interested in studying the lineability of families of nets.
	This sort of approach was initiated by J. Carmona Tapia \textit{et. al.} in \cite{CaFeFiSe} from several points of view and recently the second author in \cite{R} analyzed several ``monstrous'' families of nets and $\kappa$-sequences related to Measure Theory.
	Still about reference \cite{R}, more specifically about \cite[Section 1.2]{R}, the author showed several details regarding the study of lineability involving families of nets. Let us dive into some of these ideas since they will be relevant in our context as well.
	Given  a vector space $V$ and  a directed set $\mathcal A_0$, we assume that there exists $M_{\mathcal A_0}\subset V^{\mathcal A_0}$ satisfying a pathological property (P).
	It may be possible that we can extend $\mathcal A_0$ to a directed set $\mathcal A_1$ having greater cardinality such that we can find $M_{\mathcal A_1}\subseteq V^{\mathcal A_1}$ still satisfying property (P).
	With this mind, not only we are interested in studying the lineability properties of the family of nets satisfying property (P) indexed by an arbitrary directed set, but also we are looking for the smallest indexed set $\mathcal A$ in terms of cardinality for which there is a net indexed by $\mathcal A$ satisfying such a property and having the same lineability properties.
	Likewise, in terms of $\kappa$-sequences, we are looking for the smallest $\kappa$.

	\vspace{0.3cm} 
	\subsection{Our results} Now we are ready to describe briefly our main results. By using Fichtenholz-Kantorovich-Hausdorff theorem (see \cite{FK, Ha}), we show in Theorem \ref{thm:unbweakly} that, for an infinite-dimensional topological vector space, given a cardinal number $\kappa$ between $\aleph_0$ and $\mathfrak c$, there exists a directed set $\mathcal{A}$ of cardinality $\kappa$ such that the family of all nets indexed by $\mathcal{A}$ that are weakly convergent and unbounded is $2^{\kappa}$-lineable. 
	Next, in Proposition \ref{prop:BanachSteinhaus}, we study how big is the set of nets which do not satisfy the Banach-Steinhaus theorem. 
	More precisely, we show that, for every normed spaces $X$ and $Y$, there exists a directed set $\mathcal{A}$ with cardinality $\kappa$ between $\aleph_0$ and $\co$ such that the set of all nets of continuous linear operators $(T_a)_{a \in \mathcal{A}}$ from $X$ into $Y$ that are pointwise convergent and $\{\|T_a\|:a \in \mathcal{A}\}$ is unbounded is $2^{\kappa}$-lineable. 
	In Section~\ref{w-dense}, we study lineability properties related to the following (already mentioned) property: in every (separable) Banach space, there is a sequence $(x_n)_{n\in \mathbb N} \subseteq X$ such that $\|x_n\| \longrightarrow \infty$ as $n \rightarrow \infty$ and $\{x_n: n \in \mathbb{N} \}$ is weakly dense in $X$. We are interested in both separable and non-separable (using $\kappa$-sequences) cases. 
	To be more precise, we will study the lineability properties of the set of all $\kappa$-sequences $(x_{\alpha})_{\alpha < \kappa}$ such that $\{ \|x_{\alpha}\|: \alpha < \kappa \}$ is unbounded and $\{x_{\alpha}: \alpha < \kappa\}$ is weakly dense in $X$, where $X$ is a normed space with density character $\mbox{dens}(X)=\kappa$. Finally we study lineability properties related to a family of nets that are divergent but its series is convergent (see Theorem \ref{series}). 
	%The last section is dedicated to remarks and open problems related to the topics of the paper.
	
	%%%%%%%%%%%%%%%%%%%%%%%%%%%
	%%%%%%%%%%%%%%%%%%%%%%%%%%%
	%%%%%%%%%%%%%%%%%%%%%%%%%%%
	%%%%%%%%%%%%%%%%%%%%%%%%%%%
	%%%%%%%%%%%%%%%%%%%%%%%%%%%
	
	\section{Main Results}
	
	Since we will be dealing with different contexts as we have mentioned in the introduction, we split this section into four subsections. 
	
	\subsection{Weakly convergent and unbounded nets} 
	%In this section, 
	We will prove the following result in terms of lineability for a family of nets which are weakly convergent and unbounded. 
	Recall that a topological vector space (TVS, for short) is a vector space endowed with a topology such that vector addition and scalar multiplication are both continuous. In this case, we denote by $X^*$ its topological dual and $\sigma(X,X^*)$ the weak topology on $X$. The symbol $\mathbb{K}$ stands for the set of real or complex numbers.
	
	%\begin{theorem}\label{thm:weaklyunbounded} Let $X$ be a real or complex infinite-dimensional TVS. There exists a directed set $\A$ of cardinality $\card(X^*)$ such that the family of nets in $X$ indexed by $\A$ that are unbounded and weakly convergent is $2^\co$-lineable.
	%\end{theorem}
	
	%Let us notice that Theorem~\ref{thm:weaklyunbounded} improves \cite{CaFeFiSe}*{Theorem~2.1} by considering arbitrary infinite dimensional TVS over $\mathbb R$ or $\mathbb C$. 
	%Furthermore, and this is the most important part for us, we get $2^\co$-lineable in our case instead of only $\co$-lineable, which was the case in \cite{CaFeFiSe}*{Theorem~2.1}.
	
	\vspace{0.2cm} 
	%In order to motivate the following result, let us make some remarks regarding Theorem~\ref{thm:weaklyunbounded}.
	%Assume that $X$ is an infinite-dimnesional Banach space such that $\card(X^*) = \mathfrak c$ (for instance, this is the case when $X=c_0$).
	%Then, the directed set $\mathcal A$ that witnesses Theorem~\ref{thm:weaklyunbounded} also has cardinality $\mathfrak c$.
	%So, a very natural question is whether we can find a directed set $\mathcal A_0$ with cardinality $<\mathfrak c$ such that the family of nets indexed by $\mathcal A_0$ that are unbounded and weakly convergent still contains large dimensional vector spaces.
	%Our second result of this subsection tackles this type of problem in the context of locally convex spaces by obtaining directed sets with cardinality between $\aleph_0$ and $\mathfrak c$ such that the dimension of the desired vector could possibly be smaller than $2^{\mathfrak c}$.
	
	%Recall that a locally convex space (LCS, for short) is a TVS whose topology is defined by a family of seminorms $\mathcal P$ such that $\bigcap_{p\in \mathcal P} \{ x \colon p(x)=0 \}=\{0\}$.
	
	\begin{theorem}\label{thm:unbweakly} 
		Let $\aleph_0 \leq \kappa \leq \mathfrak c$ be a cardinal number.
		Let $X$ be a real or complex infinite-dimensional TVS. There exists a directed set $\mathcal{A}$ of cardinality $\kappa$ such that the family of nets in $X$ indexed by $\mathcal{A}$ that are unbounded and weakly convergent is $2^{\kappa}$-lineable.
	\end{theorem}
	
	Note that Theorem~\ref{thm:unbweakly} improves and generalizes \cite{CaFeFiSe}*{Theorem~2.1} by considering arbitrary infinite-dimensional TVS over $\mathbb R$ or $\mathbb C$ and decreasing the size of the index set to make it $\leq \mathfrak c$ while still having the property of being $\co$-lineability.
	
	In Remark \ref{remark1} below we provide some remarks regarding Theorem~\ref{thm:unbweakly}, which depend on the model of ZFC that we consider.
	
	\begin{remark} \label{remark1} Assume that $X$ is an infinite-dimensional TVS. On the one hand, under ZFC+CH (where CH denotes the Continuum Hypothesis), it is clear that $\aleph_0<\aleph_1=\mathfrak c$ and $2^{\aleph_0}=\mathfrak c<2^{\mathfrak c} = 2^{\aleph_1}$. Therefore, by Theorem~\ref{thm:unbweakly}, there exist directed sets $\mathcal A_0$ and $\mathcal A_1$ with cardinalities $\card(\mathcal A_0)=\aleph_0$ and $\card(\mathcal A_1)=\aleph_1$ such that the families of nets in $X$ indexed by $\mathcal A_0$ and $\mathcal A_1$ that are unbounded and weakly convergent are $\mathfrak c$-lineable and $2^{\aleph_1}$-lineable, respectively.
		So, under ZFC+CH, the size of the directed set $\mathcal A$ can affect the dimension of the desired vector space based on Theorem~\ref{thm:unbweakly}.
		
		On the other hand, under ZFC+$\neg$CH+MA (where MA denotes Martin's Axiom and $\neg$CH the negation of the CH), we have that $2^{\aleph_0}=2^{\aleph_1}$ since $\aleph_0<\aleph_1<\mathfrak c$ (see \cite[Theorem~9.5.9]{C} and \cite[Theorem~16.20]{J}). Therefore, by Theorem~\ref{thm:unbweakly} and taking $\kappa=\aleph_0$, there is a directed set $\mathcal A$ having cardinality $\aleph_0$ such that the set of nets in $X$ indexed by $\mathcal A$ being unbounded and weakly convergent is $\mathfrak c$-lineable.
		If we took $\kappa=\aleph_1$, we would increase the size of our index set but the dimension of the desired vector space would still be $\mathfrak c$; therefore, we would not obtain a larger algebraic structure even though we are increasing the size of the index set.
	\end{remark}

	In order to prove Theorem \ref{thm:unbweakly} we need to introduce some notation and remind some relevant results in our context. We start with the concept of independent families.

	\begin{definition} \label{indfamily} Let $\Gamma$ be a nonempty set. We say that a family $\mathcal{Y}$ of subsets of $\Gamma$ is {\it independent} if for any pairwise distinct sets $Y_1, \ldots, Y_n \in \mathcal{Y}$ and any $\eps_1, \ldots, \eps_n \in \{0,1\}$ we have that
		\begin{equation*}
			Y_1^{\eps_1} \cap \cdots \cap Y_n^{\eps_n} \not=\varnothing,
		\end{equation*}
		where $Y^1$ and $Y^0$ denote $Y$ and $\Gamma \setminus Y$, respectively.
	\end{definition}
	
	We will use Fichtenholz-Kantorovich-Hausdorff theorem (FKH, for short) as stated below.
	
	\begin{theorem}[Fichtenholz-Kantorovich-Hausdorff theorem] \cite{FK, Ha} \label{FKH} Let $\Gamma$ be a set of infinite cardinality $\kappa$. There is a family of independent subsets $\mathcal{Y}$ of $\Gamma$ of cardinality $2^{\kappa}$. 
	\end{theorem}
	
	It is worth mentioning the following observation below about Theorem \ref{FKH}.
	
	\begin{remark} \label{FKH-sets-are-infinite} When one applies FKH, one gets a family $\mathcal{Y}$ of $2^{\kappa}$-many subsets of nonempty sets such that $Y_1^{\eps_1} \cap \cdots \cap Y_n^{\eps_n} \not= \varnothing$ whenever $Y_1, \ldots, Y_n \in \mathcal{Y}$ and $\eps_1,\ldots,\eps_n \in \{0,1\}$. As a matter of fact, the sets $Y_1^{\eps_1} \cap \cdots \cap Y_n^{\eps_n}$ besides being nonempty are in fact infinite. The reader can go to the observation right after \cite[Definition 1.3]{FRST} for a simple proof of this fact (alternatively the reader may consider the definition of independent subsets given in the paragraph below \cite[Theorem~2.7]{Je}).
	\end{remark}

	We are now ready to provide a proof for Theorem \ref{thm:unbweakly}.

	\begin{proof}[Proof of Theorem \ref{thm:unbweakly}]
		%Take $\mathcal P$ a family of seminorms that witnesses the definition of LCS for $X$.
		Let $\mathcal A^-\subseteq (-\infty,0)$ and $\mathcal A^+ \subseteq (0,\infty)$ be such that $\card(\mathcal A^-)=\card(\mathcal A^+)=\kappa$.
		Let $\mathcal A = \mathcal A^- \cup \mathcal A^+$ endowed with the standard order of $\mathbb R$ and note that $\card(\mathcal A)= \kappa$.
		By FKH, there is a family $\mathcal K \subseteq \mathcal{P}(\mathcal A^-)$ of independent subsets of $\mathcal A^-$ such that $\card(\mathcal{K}) = 2^\kappa$. Fix $x \in X\setminus \{0\}$. For every $K \in \mathcal K$, we define 
		$$
		(x_a^K)_{a \in \mathcal A} = 
		\begin{cases}
			|a|x, & \mbox{if } a \in K, \vspace{0.2cm} \\
			\frac{1}{a} x, & \mbox{if } a \in \mathcal A^+, \vspace{0.2cm}  \\
			0, & \mbox{if } a \not\in K \cup \mathcal A^+.
		\end{cases}
		$$
		We will see that given nonzero scalars $\lambda_1, \ldots, \lambda_m \in \mathbb{K}$ and $K_1, \ldots, K_m \in \mathcal{K}$ distinct, the net $\sum_{j=1}^m \lambda_j \left(x_a^{K_j} \right)_{a\in \mathcal A}$ is unbounded, which also immediately implies the linear independence of $\left\{ \left(x_a^K\right)_{a \in \mathcal A } \colon K \in \mathcal K \right\}$. Indeed, since $\mathcal K$ is a family of independent sets, by Remark~\ref{FKH-sets-are-infinite}, there exists a sequence of distinct terms $(a_l^1)_{l\in \mathbb N} \subseteq K_1 \setminus (K_2 \cup \cdots \cup K_m)$, that is, unbounded in $\mathcal A^-$. Therefore, the set 
		\begin{equation*}
			\left\{ \sum_{i=1}^m \lambda_i x_{a_l^1}^{K_i} \colon l \in \mathbb N \right\} = \left\{ \lambda_1 \left|a_l^1\right| x \colon l \in \mathbb N \right\} 
		\end{equation*}
		is unbounded.

		Finally, it is enough to prove that each $\left(x_a^K\right)_{a \in \mathcal A }$ weakly converges to $0$. 
		Fix $x^* \in X^*$ and $K \in \mathcal K$. Since $\mathcal A$ is a linearly ordered set and $x^*(x)$ is fixed, for every $\varepsilon >0$, there is an $a_0\in \mathcal A^+ \subseteq \mathcal A$ such that for any $a\geq a_0$ we have $\frac{1}{a} \left| x^*(x) \right| < \varepsilon$.
		Hence, for any $a \geq a_0$, it yields
		\begin{equation*}
			\left|x^* \left( x_a^K \right) - x^*(0)\right| = \left| x^* \left( x_a^K \right) \right| = \frac{1}{a} \left| x^*(x) \right| < \varepsilon.
		\end{equation*}
		This finishes the proof.
	\end{proof}

	\subsection{Nets that fail the Banach-Steinhaus theorem}
	
	The idea of the construction of the vector space in the proof of Theorem~\ref{thm:unbweakly} can be carried out to show the existence of large vector spaces of nets that fail the Banach-Steinhaus theorem. 
	Although its proof uses similar ideas from the proof of Theorem~\ref{thm:unbweakly}, we provide a detailed argument for the sake of completeness. 
	
	\vspace{0.2cm} 
	Recall that the well-known Banach-Steinhaus theorem states the following (see, for instance, \cite[Chapter~3, Theorem~14.6]{Co}): let $X$ be Banach space and $Y$ be a normed space; denote by $\mathcal{L}(X,Y)$ the Banach space of all continuous linear operators from $X$ into $Y$; if a sequence $(T_n)_{n\in \mathbb N} \subseteq \mathcal{L}(X,Y)$ strongly converges pointwise, then there is a $T\in \mathcal{L}(X,Y)$ such that $(T_n)_{n\in \mathbb N}$ strongly converges pointwise to $T$ and $\{\|T_n\| \colon n\in \mathbb N \}$ is uniformly bounded. Let us recall that this theorem is a result only about sequences, not nets (for an easy example, see \cite[page 97]{Co} just after its proof as a consequence of the Principle of Uniform Boundedness).
	
	\vspace{0.2cm}
	
	In terms of lineability of the nets which do not satisfy the Banach-Steinhaus theorem, we have the following result.

	\begin{proposition}\label{prop:BanachSteinhaus}
		Let $\aleph_0 \leq \kappa \leq \mathfrak c$ be a cardinal number.
		If $X$ and $Y$ are real or complex normed spaces, then there exists a directed set $\mathcal A$ with $\card (\mathcal A)=\kappa$ such that the set of nets of continous linear operators $(T_a)_{a\in \mathcal{A}}$ in $\mathcal{L}(X,Y)$ that strongly converge pointwise and also strongly converge pointwise to an operator but $\{ \|T_a\| \colon a\in \mathcal A \}$ is not bounded is $2^{\kappa}$-lineable.
	\end{proposition}
	
	\begin{proof}
		Fix $I \in \mathcal L(X,Y) \setminus \{0\}$ (this can be done since $X\neq \{0\} \neq Y$).
		Let $\mathcal A^-\subset (-\infty,0)$ and $\mathcal A^+ \subset (0,\infty)$ be such that $\card(\mathcal A^-)=\card(\mathcal A^+)=\kappa$, and take $\mathcal A = \mathcal A^- \cup \mathcal A^+$ endowed with the standard order of $\mathbb R$.
		Note that $\card(\mathcal A)=\kappa$.
		By FKH, there exists $\mathcal K \subseteq \mathcal P(\mathcal A^-)$ a family of independent subsets of $\mathcal A^-$ having cardinality $2^{\kappa}$.
		For any $K \in \mathcal K$, define 
		$$
		\left( T_a^K \right)_{a\in \mathcal A} = \begin{cases}
			|a|I, & \text{if } a \in K,\\
			\frac{1}{a}I, & \text{if } a \in \mathcal A^+,\\
			0, & \text{otherwise}.
		\end{cases}
		$$
		Given  nonzero scalars $\lambda_1,\ldots,\lambda_m \in \mathbb{K}$  and $K_1, \ldots, K_m \in \mathcal{K}$ distinct, let us show that 
		\begin{equation*} 
			\left\{ \left\|\sum_{j=1}^m \lambda_j T_a^{K_j} \right\| \colon a \in \mathcal A \right\}
		\end{equation*} 
		is unbounded, proving also the linear independence of the subset $\left\{ \left\| \left( T_a^K \right)_{a\in \mathcal A} \right\| \colon K \in \mathcal K \right\}$.
		As $\mathcal K$ is a family of independent subsets of $\mathcal{A}^-$, we can take an unbounded sequence $(a_l^1)_{l\in \mathbb N} \subseteq K_1 \setminus (K_2 \cup \cdots \cup K_m)$, which shows that
		$$
		\left\{ \left\|\sum_{j=1}^m \lambda_j T_{a_l^1}^{K_j} \right\| \colon l \in \mathbb N \right\} = \left\{ |\lambda_1 a_l^1| \left\| I \right\| \colon l \in \mathbb N \right\}
		$$
		is unbounded.
		
		Finally, observe for any $K\in \mathcal K$ and for every $x\in X$, we have that $\left(T_a^K(x)\right)_{a\in \mathcal A}$ converges to $0$ and also $\left(T_a^K\right)_{a\in \mathcal A}$ strongly converges pointwise to the null operator.
		Thus, the vector space generated by the family of nets $\left\{ \left\| \left( T_a^K \right)_{a\in \mathcal A} \right\| \colon K \in \mathcal K \right\}$ is as needed.
	\end{proof}
	
	\subsection{Weakly dense and norm-unbounded nets} \label{w-dense} In this section, we will be interested in nets $(x_a)_{a \in \mathcal{A}}$ such that $\{x_a: a \in \mathcal{A} \}$ is weakly dense and the set $\{ \|x_a\|: a \in \mathcal{A} \}$ is unbounded. In \cite[Corollary 5]{AGM} (see also \cite{Bayart, Kadets, KLO} for more general results in this line), the authors show that, in a separable Banach space $X$, there exists a sequence $(x_n)_{n \in \mathbb{N}} \subseteq X$ such that $\{x_n: n \in \mathbb{N} \}$ is weakly dense in $X$ and $\|x_n\| \longrightarrow \infty$ as $n \rightarrow \infty$. Our aim here is to prove Theorem \ref{main7} below, which deals with a family of $\kappa$-sequences satisfying that $\{ x_{\alpha}: \alpha < \kappa \}$ is weakly dense in $X$ and also that the set $\{ \|x_{\alpha}\|: \alpha < \kappa \}$ is unbounded. We show that such a family is $\kappa^+$-lineable under some condicions on $\kappa$. 
	As an immediate consequence of Theorem~\ref{main7}, we obtain Corollary~\ref{main3} below which is related to the existence of norm divergent sequences that are weakly dense. 
	
	\begin{corollary} \label{main3} Let $X$ be real or complex separable Banach space. 
		The set of all sequences $(x_n)_{n\in \mathbb N} \subseteq X$ such that 
		\begin{itemize}
			\item[(a)] $\{ \|x_n\| : n\in \mathbb N \}$ is unbounded and 
			\item[(b)] $\overline{\{x_n: n \in \mathbb{N} \}}^{\sigma(X, X^*)} = X$
		\end{itemize}
		is $\co$-lineable.
	\end{corollary}
	
	For the proof of Theorem \ref{main7}, we need the definition of almost disjoints subsets of a regular cardinal number as well as Lemma \ref{Lemma-kappa-Jech}, which is taken from \cite{Je}. We start with almost disjointness.

	\begin{definition} \label{disjoint} Let $\kappa$ be a regular cardinal number.
		We say that $K_1,K_2 \subseteq \mathcal{P}(\kappa)$ are {\it almost disjoint} if $|K_1|=|K_2|=\kappa$ and $|K_1 \cap K_2|<\kappa$. 
	\end{definition}
	
	We will use the following lemma.
	
	\begin{lemma} \cite[Lemma 9.23 and Exercise 9.12]{Je} Let $\kappa$ be a regular cardinal.
		\begin{itemize} 
			\item[(a)]  \label{Lemma-kappa-Jech}  There exists an almost disjoint $\mathcal K \subseteq \mathcal{P}(\kappa)$ such that $|\mathcal K|=\kappa^+$.
			\item[(b)]  \label{Lemma-kappa-Jech-large} If $2^{<\kappa}=\kappa$, then there exists an almost disjoint $\mathcal K \subseteq \mathcal{P}(\kappa)$ such that $|\mathcal K|=2^\kappa$.
		\end{itemize} 
	\end{lemma}

	%In what follows, we denote the density character of the space $X$ by $\mbox{dens}(X)$. 

	\begin{theorem} \label{main7} Let $X$ be a real or complex normed space with $\mbox{dens}(X) = \kappa\geq \aleph_0$, where $\kappa$ is a regular cardinal. 
		Denote by $\textup{UWD}_\kappa$ the set of all $\kappa$-sequences $(x_{\alpha})_{\alpha < \kappa}$ such that
		\begin{itemize}
			\item[\textup{(i)}] $\{ \|x_{\alpha}\|: \alpha < \kappa \}$ is unbounded and 
			\item[\textup{(ii)}] $\overline{ \{x_{\alpha}: \alpha < \kappa \}}^{\sigma(X, X^*)} = X$.
		\end{itemize}
		Then, the set $\textup{UWD}_\kappa$ is $\kappa^+$-lineable.
		Moreover, if $2^{<\kappa}=\kappa$, then $\textup{UWD}_\kappa$ is $2^{\kappa}$-lineable.
	\end{theorem}
	
	Before we get into the proof of Theorem \ref{main7}, let us provide a comment about properties (i) and (ii) above. 
	Suppose that the net $(x_{a})_{a\in \mathcal A}$ satisfies that $\{ x_a : a\in \mathcal A \} = X$ or $\{ x_a : a\in \mathcal A \}$ is norm-dense in $X$. Then, clearly we have that $\{ \|x_a\| : a\in \mathcal A \}$ is unbounded and $\{x_a : a\in \mathcal A \}$ is weakly dense (i.e., conditions (i) and (ii) are both satisfied for the net $(x_{a})_{a\in \mathcal A}$). 
	Moreover, if $\mbox{dens}(X) = \kappa$, we can assume that $|\mathcal A|=\kappa$ and consider a bijection $f: \mathcal{A} \longrightarrow \kappa$. 
	Since (i) and (ii) are satisfied for $(x_a)_{a \in \mathcal{A}}$, conditions (i) and (ii) will be also satisfied for a $\kappa$-sequence thanks to the bijection $f$. 
	This guarantees the existence of $\kappa$-sequences satisfying such properties.
	
	\begin{proof}[Proof of Theorem \ref{main7}] 
		We start the proof by defining a family $\mathcal K$ of almost disjoint subsets (see Definition \ref{disjoint} above) of $\kappa$ that witnesses either Lemma~\ref{Lemma-kappa-Jech}.(a) or Lemma~\ref{Lemma-kappa-Jech-large}.(b). 
		For each $K \in \mathcal{K}$, we can consider an increasing bijection $\varphi_K: K \longrightarrow \kappa$. 
		Now, let $(x_{\alpha})_{\alpha < \kappa}$ be a $\kappa$-sequence satisfying properties (i) and (ii). 
		For a fixed $K \in \mathcal{K}$ and every $\alpha < \kappa$, we define 
		\begin{equation} \label{newelement}
			x_{\alpha}^K := \begin{cases}
				x_{\varphi_K(\alpha)}, & \mbox{whenever } \alpha \in K,\\
				0, & \mbox{otherwise}. 
			\end{cases}
		\end{equation}
		Write $K = \{\alpha_{\beta}: \beta < \kappa \}$, where $\beta_1 < \beta_2<\kappa$ implies that $\alpha_{\beta_1} < \alpha_{\beta_2}$. 
		Since each $K \in \mathcal{K}$ is such that $|K| = \kappa$, we have that the set $\{ \|x_{\alpha}^K\|: \alpha < \kappa \}$ is also unbounded. Moreover, since for every $K \in \mathcal{K}$ we have that 
		\begin{equation*}
			\{ x_{\alpha}: \alpha < \kappa \} \subseteq \{ x_{\alpha}^K: \alpha < \kappa \},
		\end{equation*}
		then we obtain
		\begin{equation*}
			X = \overline{ \{ x_{\alpha}: \alpha < \kappa \}}^{\sigma (X, X^*)} \subseteq \overline{ \{ x_{\alpha}^K: \alpha < \kappa \}}^{\sigma (X,X^*)} \subseteq X.
		\end{equation*}
		This shows that $\overline{ \{ x_{\alpha}^K: \alpha < \kappa \}}^{\sigma (X,X^*)} = X$ for every $K \in \mathcal{K}$.
		
		\vspace{0.3cm} 
		
		Now, let $\lambda_1, \ldots, \lambda_n$ be nonzero scalars and $K_1, \ldots, K_n \in \mathcal{K}$ distinct. 
		Since $\card(K_1) = \kappa$ and $\card(K_1 \cap K_j) < \kappa$ for every $j=2,\ldots, n$, we have that the cardinality of the set
		\begin{equation*}
			\mathcal{I}_1 := K_1 \setminus \left( \bigcup_{j=2}^{n} (K_1 \cap K_j) \right)  
		\end{equation*}
		is also $\kappa$, that is, $\card(\mathcal{I}_1) = \kappa$. 
		Consider $(x_{\alpha_{\beta}})_{\beta < \kappa}$ to be a $\kappa$-subsequence of $(x_{\alpha})_{\alpha < \kappa}$ indexed by $\mathcal{I}_1$. Then, by the election of the index set $\mathcal{I}_1$ and by using (\ref{newelement}), we get that 
		\begin{equation*}
			\left\{ \left\| \sum_{j=1}^n \lambda_j x_{\alpha_{\beta}}^{K_j} \right\|: \beta < \kappa \right\} 
			%= \left\{ \lambda_1 x_{\alpha_{\beta}}^{K_1}: \beta < \kappa \right\} 
			= \left\{ |\lambda_1| \left\| x_{\varphi_{K_1}(\beta)} \right\|: \beta < \kappa \right\}.
		\end{equation*}
		This argument shows at once that $(x_{\alpha}^{K_1})_{\alpha < \kappa}, \ldots, (x_{\alpha}^{K_n})_{\alpha < \kappa}$ are linearly independent and also that the set 
		\begin{equation*}
			\left\{ \left\| \sum_{j=1}^n \lambda_j x_{\alpha}^K \right\|: \alpha < \kappa \right\} 
		\end{equation*}
		is unbounded in $X$. 
		It remains to prove that 
		\begin{equation} \label{1}
			\overline{ \left\{ \sum_{j=1}^n \lambda_j x_{\alpha}^{K_j}: \alpha < \kappa \right\} }^{\sigma (X,X^*)} = X.
		\end{equation}
		Before doing that, let us observe the following. 
		Since $\card \left( \bigcup_{j=2}^n (K_1 \cap K_j)\right) < \kappa$ and $\kappa$ is regular it yields
		\begin{equation*}
			\overline{\alpha} := 
			%\sup \{ \alpha_j: 2 \leq j \leq k \} = 
			\sup \left\{ \alpha \in \bigcup_{j=2}^n (K_1 \cap K_j ) \right\} < \kappa
		\end{equation*}
		by \cite[Lemma~3.9]{Je}.
		Why we are considering such an $\overline{\alpha}$ will be clear below. 
		In order to prove (\ref{1}), let us fix $x \in X$. Since (ii) holds true, there exists $(z_{\beta})_{\beta < \kappa} \subseteq \{ x_\alpha : \alpha < \kappa \}$ such that $z_{\beta} \stackrel{w}{\longrightarrow} \frac{1}{\lambda_1} x$. Considering once again the index set $\mathcal{I}_1$ and $\overline{\alpha}$, there exists $\overline{\alpha} \leq \alpha_0<\kappa$ such that for every $\alpha_0\leq \beta<\kappa$, we have that 
		\begin{equation*}
			\sum_{j=1}^n \lambda_j x_{\beta}^{K_j} = \lambda_1 x_{\varphi_{K_1}(\beta)}.
		\end{equation*}
		Since $\varphi_{K_1}$ is an increasing bijection and $\card(\mathcal{I}_1) = \kappa$, we have that 
		\begin{equation*}
			\left\{ \lambda_1 x_{\varphi_{K_1}(\beta)} : \alpha_0 \leq \beta < \kappa \right\} = \{ \lambda_1 x_{\alpha}: \alpha < \kappa \} \setminus \mathcal{F},
		\end{equation*}
		where $\mathcal{F}$ is a set with elements of the $\kappa$-sequence $(\lambda_1 x_{\alpha})_{\alpha < \kappa}$ such that its cardinality is $< \kappa$. 
		Therefore, there exists a $\kappa$-sequence $(\tilde{z}_{\beta})_{\beta < \kappa} \subseteq \{ \sum_{j=1}^n \lambda_j x_{\alpha}^{K_j}: \alpha < \kappa \}$ such that 
		\begin{equation*}
			\tilde{z}_{\beta} \stackrel{w}{\longrightarrow} \lambda_1 \cdot \frac{1}{\lambda_1}x = x 
		\end{equation*}
		and this proves (\ref{1}) as desired. 
	\end{proof}

	\subsection{Convergent series with associated divergent nets}

	Let $X$ be a normed space and $\mathcal I$ an infinite set.
	We can give meaning to the convergence of the (possibly) ``uncountable sum'' in $X$, denoted by $\sum_{i\in I} x_i$, where each $x_i$ belongs to $X$, as follows: consider $\mathcal F$ to be the set of all finite subsets of $I$ endowed with the inclusion $\subseteq$. Bearing this in mind, we have that $\mathcal F$ is a directed set.
	Now, for every $F\in \mathcal F$, we define
	$$
	x_F := \sum_{i\in F} x_i.
	$$
	Each $x_F$ is then a well-defined vector of $X$ (since $F$ is finite) and $(x_F)_{F\in \mathcal F}$ is a net. In the same line, we have the following definition. Given $x_i \in X$ for all $i\in I$, we say that $\sum_{i\in I} x_i$ converges to $x\in X$ whenever $\lim_{F\in \mathcal F} x_F = x$. Recall that in Hilbert spaces, the latter definition can be used to obtain some relevant characterization in the non-separable case (see, for instance, \cite[Chapter~1, Theorem~4.13]{Co}).

	Our next result is motivated by the following fact. If $X=H$ is a Hilbert space and $I=\mathbb N$, it is known that if $\lim_{F\in \mathcal F} h_F = h \in \mathcal H$, then $\sum_{n=1}^\infty h_n = h$, but the converse is not true in general (it is important to mention that if $\sum_{n=1}^\infty h_n$ is absolutely convergent, then the converse implication is satisfied; see, for instance, \cite[Chapter~1, Section~4, Exercises~10 and~11]{Co}). The following result shows that we can find (in terms of lineability) large sets of sequences $(x_n)_{n\in \mathbb N}$ in a normed space such that the series $\sum_{n=1}^\infty x_n$ is convergent but the net $(x_F)_{F\in \mathcal F}$ diverges.
	In what follows, we denote by $\mbox{CS}(\mathbb K) \subseteq \mathbb K^{\mathbb N}$ the set of all sequences $(k_n)_{n\in \mathbb N}$ such that the series with general terms $(k_n)_{n\in \mathbb N}$ is conditionally convergent. 
	
	\begin{theorem} \label{series} Let $X$ be a normed space defined over $\mathbb K \in \{ \mathbb R,\mathbb C \}$ and $\mathcal F$ the family of finite subsets of $\mathbb N$ endowed with the order $\subseteq$. The set of all sequences $(x_n)_{n\in \mathbb N} \subseteq  X$ such that $\sum_{n=1}^\infty x_n$ is convergent and $(x_F)_{F\in \mathcal F}$ diverges is $\mathfrak c$-lineable.
	\end{theorem}
	
	\begin{proof}
		Since a series $\sum a_n$ of complex numbers is conditionally convergent if and only if $\sum \mbox{Re}(a_n)$ or $\sum \mbox{Im}(a_n)$ is conditionally convergent, we restrict ourselves to the real case, that is, $\mathbb{K} = \mathbb{R}$. Fix $x\in X\setminus \{0\}$.
		Let $V_1\subseteq (\text{CS}(\mathbb K) \cup \{0\})$ be a vector subspace of dimension $\mathfrak c$ (that we might consider thanks to \cite[Theorem 2.1]{APS}, which holds for both real and complex cases).
		Then, the set 
		\begin{equation*} 
			xV_1 := \{ (k_n x)_{n\in \mathbb N} \colon (k_n)_{n\in \mathbb N} \in V_1 \} \subseteq X^{\mathbb N}
		\end{equation*}
		yields the desired result.
		Indeed, it is easy see that $xV_1$ is a vector subspace of $X^{\mathbb N}$ of dimension $\mathfrak c$ such that $\sum_{n=1}^\infty k_n x = x \sum_{n=1}^\infty k_n$ converges.
		Recall that given a conditionally convergent series $\sum_{n=1}^\infty k_n$, the series of positive terms $\sum_n k_n^+$ diverges. 
		Now fix $M>0$ and $F \in \mathcal F$ to be arbitrary. 
		Let $F^+_M \in \mathcal F$ be such that $F_M^+ \supseteq F$, $k_n=k_n^+$ provided that $n\in F^+_M \setminus F$, and 
		$$
		\sum_{n\in F_M^+} k_n = \sum_{n\in F} k_n + \sum_{n\in F_M^+\setminus F} k_n > M.
		$$
		Therefore,
		$$
		\left\| \sum_{n\in F_M^+} k_n x \right\| = \left\| x \sum_{n\in F_M^+} k_n \right\| = \left| \sum_{n\in F^+_M} k_n \right| \|x\| > M\|x\|.
		$$
	\end{proof}

	\noindent 
	\textbf{Acknowledgements:} The authors are thankful to Manuel Maestre (Universidad de Valencia), Vladimir Kadets (Holon Institute of Technology) and Juan B. Seoane-Sepúlveda (Universidad Complutense de Madrid) for useful conversations about different parts of the paper during the writing procedure.
	This paper was partially completed at the Universidad de Valencia (UV) during a research stay of the second author. 
	He expresses his gratitude for the support and warm reception received by the Departamento de Análisis Matemático at UV.
	
	\vspace{0.2cm} 
	\noindent 
	\textbf{Funding information}: In what follows, we provide the list of funding that the authors have received during the period this paper was being prepared.

	\vspace{0.2cm}
	\noindent
	{\it Sheldon Dantas} was supported by 
	\begin{itemize} 
		\item[(1)] The Spanish AEI Project PID2019 - 106529GB - I00 / AEI / 10.13039/501100011033.
		\item[(2)]  Generalitat Valenciana project CIGE/2022/97.
		\item[(3)] The grant PID2021-122126NB-C33 funded by MCIN/AEI/10.13039/501100011033 and by “ERDF A way of making Europe”. 
	\end{itemize} 
	\vspace{0.2cm}
	\noindent
	{\it Daniel L. Rodríguez-Vidanes} was supported by 
	\begin{itemize}
		\item[(1)] The grant PGC2018-097286-B-I00. 
		
		\item[(2)] The Spanish Ministry of Science, Innovation and Universities and the European Social Fund through a “Contrato Predoctoral para la Formación de Doctores, 2019” (PRE2019-089135).
	\end{itemize}

\end{document}